\documentclass{article}
\usepackage{amsmath}
\usepackage{amssymb}
\usepackage{amsfonts}
\usepackage{amsthm}
\usepackage{hyperref}
\usepackage{color}
\usepackage{graphicx} % Required for inserting images
\usepackage{algorithm}
\usepackage{algpseudocode}
\usepackage{booktabs}
\usepackage{xr-hyper}
\usepackage{algorithm}
\usepackage{algpseudocode}
\usepackage{subfiles}
\usepackage{booktabs,multirow}% provides \State, \Procedure, etc.
\def\E{\mathbb{E}}

\usepackage[style=authoryear,maxcitenames=1,maxbibnames=99]{biblatex}
\usepackage{geometry}
\newtheorem{theorem}{Theorem}[section]

\newtheorem{proposition}{Proposition}[section]

\newtheorem{lemma}{Lemma}[section]

\newtheorem{definition}{Definition}[section]

\newtheorem{remark}{Remark}[section]

\newtheorem{assumption}{Assumption}[section]

\title{Double/debiased machine learning of quantile treatment effects on long-term outcomes in clinical trials}
	\author{Ziyang Liu$^1$, Niwen Zhou$^1$, Peng Wu$^{2}$,   and Xu Guo$^3$
	\\~
  {\small \it $^{1}$ Faculty of Arts and Sciences, Beijing Normal University, Zhuhai, China}\\
   {\small \it $^{2}$ 
School of Mathematics and Statistics, Beijing Technology and Business University, Beijing, China}\\
  {\small \it $^{3}$ School of Statistics, Beijing Normal University, Beijing, China}  
	}
	\date{}
\date{ }

\begin{document}

\maketitle

\begin{abstract}
Long-term outcomes are often unavailable in randomized clinical trials, although short-term surrogate outcomes are commonly observed. External observational data may contain the long-term outcome, but causal comparisons based on such data alone are vulnerable to confounding. Existing surrogate-based data integration methods for long-term outcomes have focused primarily on average treatment effects.
We study estimation of quantile treatment effects for long-term outcomes in the trial population by combining randomized trial data with external observational data. Under treatment randomization, positivity, and a surrogate-based transportability assumption, we establish identification and develop a doubly robust estimator for inference. The estimator accommodates flexible machine learning methods for nuisance estimation, remains consistent if either the score-related or outcome regression-related nuisance functions are consistently estimated, and is asymptotically normal under regularity conditions.
Simulation and real-data results demonstrate that the proposed method performs well in finite samples and can reveal heterogeneous long-term treatment effects across quantiles.

\end{abstract}

\section{Introduction}

Randomized controlled trials (RCTs) are widely regarded as the gold standard for causal inference in clinical research because treatment assignment is randomized by design. In many trials, however, the primary outcome of scientific interest is long-term and cannot be fully observed within the trial itself. Extended follow-up is often costly and logistically difficult, and outcome ascertainment may require linkage to external data sources such as disease registries, electronic health records, or administrative claims. As a result, clinical trials frequently measure short-term biomarkers or intermediate outcomes, while the long-term clinically relevant outcome remains unavailable in the experimental sample. By contrast, external observational data sources may contain long-term outcome information at a much larger scale, but they do not generally provide randomized treatment assignment and are therefore vulnerable to unmeasured confounding. This creates a central challenge for long-term causal inference: the trial provides credible treatment variation but limited long-term outcome information, whereas observational data contain the long-term outcome but may not support unbiased causal comparisons on their own. Throughout, our target is the long-term causal effect in the trial population, with the observational sample serving only as an auxiliary source of information on the primary outcome.

A natural strategy is therefore to combine these two data sources in a way that exploits their complementary strengths. In particular, recent work has studied how short-term surrogate outcomes measured in the trial can be used to bridge the experimental sample to auxiliary observational data when the primary long-term outcome is unavailable or only partially observed in the trial. Under suitable assumptions linking the long-term outcome to the observed surrogate, baseline covariates, and treatment across samples, these approaches use randomized treatment variation from the trial together with long-term outcome information from external data to recover causal effects for the trial population. This surrogate-based data integration framework has become an increasingly important approach to long-term causal inference when direct experimental follow-up is costly or infeasible; see, for example, \textcite{chen2023semiparametric,athey2025surrogate,kallus2025role,hu2025identification,imbens2025long}.

Despite this progress, the existing literature has focused primarily on mean effects, especially the average treatment effect (ATE). In many applications, however, the mean provides only a limited summary of treatment impact, particularly when long-term outcomes are skewed or heterogeneous across the distribution. Quantile treatment effects (QTE) are therefore of direct interest because they capture how treatment shifts the distribution of long-term outcomes beyond the mean. Moreover, in the surrogate-based data integration setting, QTE are not a routine extension of mean effects, because quantiles depend on the full outcome distribution and therefore create distinct challenges for both identification and estimation.

This distinction is also practically important in our motivating application to IgA nephropathy (IgAN). Under the same clinical setting, \textcite{liu2019effects} showed that hydroxychloroquine (HCQ) combined with renin--angiotensin--aldosterone system (RAAS) inhibition improves the surrogate outcome, and \textcite{hu2025identification} further found a favorable average treatment effect on the long-term renal outcome relative to standard RAAS therapy alone. While such findings are informative, they do not reveal whether the treatment benefit is homogeneous across patients with different levels of long-term renal deterioration. In contrast, our analysis of the long-term quantile treatment effect shows that HCQ combined with RAAS therapy is associated with lower long-term renal deterioration over most parts of the outcome distribution, and that this benefit becomes more pronounced at higher quantiles. These findings therefore provide evidence of substantial treatment effect heterogeneity. In particular, they suggest that HCQ combined with RAAS therapy may offer stronger protection for patients who would otherwise experience more severe renal decline, a feature that cannot be fully captured by the average treatment effect alone.

Quantile treatment effects have been extensively studied in standard causal settings where the outcome is observed in the target sample. Existing work has developed semiparametrically efficient, doubly robust, and machine-learning-based methods for QTE estimation; see, for example, \textcite{firpo2007efficient,frolich2013unconditional,diaz2017efficient,kallus2024localized,cheng2024doubly}. However, these methods do not directly apply to our setting, in which the long-term outcome is unavailable in the experimental sample and must instead be recovered through surrogate-based integration of experimental and observational data. To the best of our knowledge, the identification and doubly robust estimation of long-term QTE in this surrogate-based data integration setting have not been studied.

We make four main contributions. First, we establish identification results for the long-term QTE under surrogate-based assumptions linking the experimental and observational samples. Second, we derive the efficient influence function and the corresponding semiparametric efficiency bound for the target parameter in this nonparametric model. Third, we develop a doubly robust estimator that permits flexible estimation of nuisance functions and is compatible with modern machine learning methods. Fourth, we establish asymptotic theory and valid large-sample inference for the proposed estimator. Taken together, our results extend surrogate-based long-term causal inference from mean effects to quantile effects and connect that literature with semiparametric methods for distributional causal inference. We also illustrate the practical value of the proposed framework through a real-data application to an IgA nephropathy study, where randomized trial data are combined with external observational data to assess treatment effects on long-term renal outcomes.

The rest of the paper is organized as follows. Section~\ref{sec:Preliminaries} introduces the setup, notation, and identifying assumptions. Section~\ref{sec:Semiparametric Efficiency} derives the efficient influence function and the semiparametric efficiency bound for the long-term quantile treatment effects. Section~\ref{sec:Estimation} presents the doubly robust estimation procedure and develops its asymptotic properties. Section~\ref{sec: Simulation} reports simulation results. Section~\ref{sec:Application} presents a real-data application. Section~\ref{sec:Conclusion} concludes.

\section{Preliminaries}\label{sec:Preliminaries}

Let $T \in \{0,1\}$ denote a binary treatment indicator, where $T=1$ represents receiving treatment and $T=0$ represents the control condition. Let $X \in \mathcal{X} \subseteq \mathbb{R}^{d_x}$ denote baseline covariates measured before treatment assignment, and $Y \in \mathbb{R}$ represent the outcome of primary interest. Following the Neyman-Rubin potential outcomes framework, we assume two potential outcomes exist, denoted by $Y(1)$ and $Y(0)$, which correspond to the outcomes under treatment and control conditions, respectively\parencite{rubin1974estimating,splawa1990application}. By the standard consistency assumption, the observed outcome $Y=Y(T)=T\cdot Y(1)+(1-T)\cdot Y(0)$.

In this paper, we consider the setting where the primary outcome $Y$ cannot be observed for units in the RCT. However, a surrogate vector $S \in \mathcal{S} \subseteq \mathbb{R}^{d_s}$ are observed for all units in both data sources. Analogously, by the standard consistency assumption, the observed outcome $S=S(T)=T\cdot S(1)+(1-T)\cdot S(0)$.

Let \(G \in \{0,1\}\) be an indicator of sample membership, where \(G=1\) denotes the experimental sample and \(G=0\) denotes the observational sample. For each unit, the observed data vector is \(W=(G,T,X,S,(1-G)Y)\), so that the primary outcome \(Y\) is observed only when \(G=0\). Equivalently, we observe an observational sample \(\{(Y_i,T_i,X_i,S_i,G_i=0): i\in\mathcal{I}^o\}\) and an experimental sample \(\{(Y_i=\mathrm{NA},T_i,X_i,S_i,G_i=1): i\in\mathcal{I}^e\}\), where \(\mathrm{NA}\) denotes a missing value. The index sets \(\mathcal{I}^o\) and \(\mathcal{I}^e\) correspond to the observational and experimental samples, respectively. We denote the sample sizes by \(N_0=|\mathcal{I}^o|\), \(N_1=|\mathcal{I}^e|\), and \(N=N_0+N_1\).

The parameter of interest is the QTE, defined as:

\begin{align*}
\Delta_\tau = q^1_{\tau} - q^0_{\tau},
\end{align*}
where $q^t_{\tau}$ is the $\tau$-quantile of the potential outcome $Y(t)$ in the experimental data, satisfying the condition:

\begin{align*}
\tau = \Pr[Y(t) \leq q^t_{\tau} \mid G=1], \quad t = 0,1.
\end{align*}

Equivalently, the $\tau$-quantile $q^t_{\tau}$ can be expressed as

\begin{align*}
q^t_{\tau} = \arg\min_{q} \mathbb{E}[\rho_\tau(Y(t)-q) \mid G=1], \quad t=0,1,
\end{align*}
where $\rho_\tau(u) = u(\tau - \mathbb{I}(u < 0))$ denotes the quantile loss function.

To identify $\Delta_\tau$, the following assumptions are required.

\begin{assumption}[Unconfoundedness of experimetal sample]\label{unconfoundedness} For $t=0,1$,

\begin{align*}
    T{\bot}(Y(t),S(t))|X,G=1. 
\end{align*} 

\end{assumption}

Assumption~\ref{unconfoundedness} states that treatment assignment is randomized conditional on baseline covariates in the experimental sample, which is standard in the potential outcomes framework for randomized studies \parencite{rubin1974estimating,HernanRobins2020}. In a completely randomized trial, this assumption holds by design. We formulate it conditional on \(X\) to accommodate more general randomization schemes, such as stratified or covariate-adaptive randomization, and to align with the subsequent estimation framework . No analogous unconfoundedness assumption is imposed on the observational sample, because our target parameter pertains to the experimental population and the observational data are used only to supplement information on the unobserved long-term outcome, rather than to identify treatment effects within the observational sample itself.

\begin{assumption}[Positivity]\label{Positivity}
There exists $\varepsilon \in (0,1/2)$ such that almost surely,
\begin{enumerate}
    \item $\varepsilon \leq \Pr(T=1 \mid S,X,G=0) \leq 1-\varepsilon$;
    \item $\varepsilon \leq \Pr(T=1 \mid X,G=1) \leq 1-\varepsilon$;
    \item $\varepsilon \leq \Pr(G=1 \mid S,X,T=t) \leq 1-\varepsilon$, \quad t=0,1.
\end{enumerate}
\end{assumption}

Assumption~\ref{Positivity} imposes the overlap conditions required for identification, which are standard in the causal inference and transportability literature  \parencite{rosenbaum1983central,hu2025identification,kallus2025role}. 
Part~1 ensures treatment overlap in the observational sample conditional on $(S,X)$, so that both treatment arms are represented for the relevant values of the covariates and surrogate. 
Part~2 ensures treatment overlap in the experimental sample conditional on $X$, which is needed to define the target causal effect in the trial population. 
Part~3 ensures overlap between the observational and experimental samples conditional on $(S,X,T)$, so that information learned from the observational sample can be validly transported to the experimental sample. 
Moreover, the marginal bound $\varepsilon \leq \Pr(G=1) \leq 1-\varepsilon$ follows directly from Part~3.

\begin{assumption}[Missing at random]\label{missing at random}For $t=0,1$,

\begin{align*}
    G{\bot}Y(t)|S(t),X,T.
\end{align*}

\end{assumption}
Assumption~\ref{missing at random} is also adopted in \textcite{kallus2025role}. This assumption
 posits that the primary outcome is missing at random conditionally. That is, the missingness indicator $G$ depends only on observed variables, including pretreatment covariates $X$, the surrogates $S$, and the treatment $T$. This assumption ensures that, after conditioning on the observed variables, the distribution of the primary outcome in the RCT and the observational data is comparable. Consequently, information from the observational data can be validly used to infer properties of the missing primary outcome in the RCT data.

\begin{assumption}[Quantile uniqueness]\label{Monotonicity}

There are nonempty sets $\mathcal{Y}_1$ and $\mathcal{Y}_0$, such that $\mathcal{Y}_t = \{\tau \in (0,1)$: $ \Pr[Y(t) \leq q^t_{\tau} - c\mid G=1] < \Pr[Y(t) \leq q^t_{\tau} + c\mid G=1], \forall c \in \mathbb{R}, c > 0\}$.

\end{assumption}

Assumption~\ref{Monotonicity} imposes strict local monotonicity of the marginal distribution function conditional on \(G=1\) around \(q^t_{\tau}\), thereby ensuring that the \(\tau\)-th marginal quantile \(q^t_{\tau}\) is uniquely defined. Such a quantile uniqueness condition is standard in the quantile treatment effect literature; see, for example, \textcite{firpo2007efficient}.

Below we show the identification of $\Delta_\tau$.

\begin{proposition}\label{identification}
Under Assumptions~\ref{unconfoundedness}--\ref{Monotonicity}, the long-term quantile treatment effect
\(\Delta_\tau=q_\tau^1-q_\tau^0\) can be identifiable as: 
\[
q_\tau^t
=
\arg\min_q\;
\mathbb{E}\left[
\mathbb{E}\left[
\mathbb{E}\left[\rho_\tau(Y-q)\mid S,X,T=t,G=0\right]
\mid X,T=t,G=1
\right]
\mid G=1
\right],
\] for $t=0, 1$.
\end{proposition}

\section{Semiparametric Efficiency}\label{sec:Semiparametric Efficiency}

In this section, we derive the efficient influence function and the semiparametric efficiency bound for the long-term quantile treatment effect \(\Delta_\tau\). 
%\subsection{Nuisance functions}\label{subsec:Nuisance functions}
The efficient influence function depends on a collection of nuisance functions, thus we first introduce the following definitions. 

The first class of nuisance functions consists of outcome bridge functions that connect the observational sample, in which the long-term outcome is observed, to the experimental population, for which the target parameter is defined. For each \(t=0,1\), let
\[
F_t(q\mid S,X)
:=
\Pr(Y\le q \mid S,X,T=t,G=0)
\]
denote the conditional distribution function of the primary outcome in the observational sample, and define the transported bridge function
\[
m_t(X,q)
:=
\mathbb{E}\!\left[F_t(q\mid S,X)\mid X,T=t,G=1\right].
\]
Here, \(F_t(q\mid S,X)\) captures the conditional distribution of the long-term outcome in the observational sample, while \(m_t(q,X)\) transports this distributional information to the experimental population by averaging over the surrogate distribution in the experimental sample.

For later use, we also define
\[
f_t(q\mid S,X)
:=
\frac{\partial}{\partial q}F_t(q\mid S,X),
\]
whenever the derivative exists.

The second class of nuisance functions consists of the treatment propensity score in the experimental sample and the sample membership score. Specifically, define
\[
e(X):=\Pr(T=1\mid X,G=1),
\]
and, for each \(t=0,1\),
\[
g_t(S,X):=\Pr(G=1\mid S,X,T=t).
\]
For notational convenience, we define
\[
\alpha_t(S,X):=\frac{g_t(S,X)}{[1-g_t(S,X)][t\,e(X)+(1-t)\big(1-e(X)\big)]},
\qquad t=0,1.
\]

For brevity, let \(\eta := (F_t, m_t, e, g_t )\) denote the collection of nuisance functions appearing in the efficient influence function. We defer the discussion of how these nuisance functions are estimated in practice to Remark~\ref{rem: estimating detail}.

%\subsection{Influence function and semiparametric efficiency bound}\label{subsec: Influence function and semiparametric efficiency bound}

Having defined the necessary nuisance functions, we are now ready to characterize the EIF for the quantile treatment effect \(\Delta_\tau\), following the approach developed in \textcite{hines2022demystifying}. The EIF consists of two components corresponding to the treated and control potential outcomes, respectively. Each component combines information from both the experimental sample and the observational sample through inverse probability weighting and transported bridge functions.

\begin{theorem}[EIF of \(\Delta_\tau\)]\label{eif}
Let \(\nu=\Pr(G=1)\). Under Assumption~\ref{unconfoundedness}--\ref{Monotonicity}, the efficient influence function of \(\Delta_\tau=q_\tau^1-q_\tau^0\) is
\[
\phi=\phi_1-\phi_0,
\]
where
\begin{align*}
\phi_t&=J_t^{-1}\psi_t(W;q^t_\tau,\eta)
\end{align*}
with
\begin{align*}
\psi_1(W;q^1_\tau,\eta)
&:=
\frac{G}{\nu}
\left[
\frac{T\left\{F_1(q^1_{\tau}\mid S,X)-m_1(X,q^1_{\tau})\right\}}{e(X)}
+m_1(X,q^1_{\tau})-\tau
\right]\\
&\quad+
\frac{T(1-G)}{\nu}\alpha_1(S,X)\left\{\mathbb{I}(Y\le q^1_{\tau})-F_1(q^1_{\tau}\mid S,X)\right\},
\end{align*}
and
\begin{align*}
\psi_0(W;q^0_\tau,\eta)
&:=
\frac{G}{\nu}
\left[
\frac{(1-T)\left\{F_0(q^0_{\tau}\mid S,X)-m_0(X,q^0_{\tau})\right\}}{1-e(X)}
+m_0(X,q^0_{\tau})-\tau
\right]\\
&\quad+
\frac{(1-T)(1-G)}{\nu}\alpha_0(S,X)\left\{\mathbb{I}(Y\le q^0_{\tau})-F_0(q^0_{\tau}\mid S,X)\right\}.
\end{align*}
Moreover,
\begin{align*}
J_t
&:=
\mathbb{E}\!\left[\mathbb{E}\!\left[f_t(q^t_{\tau}\mid S,X)\mid X,T=t,G=1\right]\mid G=1\right].
\end{align*}
Hence, the semiparametric efficiency bound is \(\mathbb{E}[\phi^2]\).
\end{theorem}

Theorem~\ref{eif} characterizes the fundamental efficiency limit for estimating the long-term quantile treatment effect \(\Delta_\tau\) under Assumptions~\ref{unconfoundedness}--\ref{Monotonicity}. In particular, for any regular estimator of \(\Delta_\tau\), the asymptotic variance cannot be smaller than \(\mathbb{E}[\phi^2]\). In the next subsection, \(\phi\) will serve as the basis for constructing our estimator.

\section{Estimation and inference}\label{sec:Estimation}

\subsection{K-fold cross-fitting estimator}

Note that the quantile treatment effect is $\Delta_\tau = q^1_{\tau}-q^0_{\tau}$. 
Building on the identification result, we propose a doubly robust estimator based on $K$-fold sample splitting to consistently estimate $\Delta_\tau$. 
For each fold $k$, let $\hat e_{-k}(X)$, \(\hat{\alpha}_{1,-k}(S,X)\) and $\hat Q_{\tau,-k}^T(\cdot)$ denote nuisance estimators trained on the complementary sample $\{1,\dots,n\}\setminus\mathcal{I}_k$. 
Let

    \begin{align*}
&\hat{Q}_{\tau,-k}^{t}(X,q)=\hat{\E}_{-k}\!\Bigl[
\hat{\E}_{-k}\bigl[\,\rho_{\tau}(Y-q)\,\big|\,S,X,T=t,G=0\bigr]
\bigm| X,T=t,G=1
\Bigr].\\
&\hat{Q}_{\tau,-k}^t(S,X,q) = \hat{\mathbb{E}}_{-k}[\rho_\tau(Y - q) \mid S,X,T=t,G=0].
    \end{align*}

\begin{definition}[doubly robust estimator]\label{def:dr estimator}
For each \(t=0,1\), define
\begin{align*}
\hat q_{dr}^t
=
\arg\min_q\;
\frac{1}{K}\sum_{k=1}^K \frac{1}{|\mathcal I_k|}\sum_{i\in\mathcal I_k}
\Psi_t(W_i;q,\hat\eta_{-k}),
\label{eq:dr-estimator}
\end{align*}
where
\begin{align*}
\Psi_1(W_i;q,\hat\eta_{-k})
&=
\frac{G_i}{\nu}
\left[
\frac{T_i\bigl\{\hat Q^{1}_{\tau,-k}(S_i,X_i,q)-\hat Q^{1}_{\tau,-k}(X_i,q)\bigr\}}{\hat e_{-k}(X_i)}
+\hat Q^{1}_{\tau,-k}(X_i,q)
\right] \\
&\quad
+\frac{T_i(1-G_i)}{\nu}\hat\alpha_{1,-k}(S_i,X_i)
\left\{
\rho_\tau(Y_i-q)-\hat Q^{1}_{\tau,-k}(S_i,X_i,q)
\right\},
\\[1ex]
\Psi_0(W_i;q,\hat\eta_{-k})
&=
\frac{G_i}{\nu}
\left[
\frac{(1-T_i)\bigl\{\hat Q^{0}_{\tau,-k}(S_i,X_i,q)-\hat Q^{0}_{\tau,-k}(X_i,q)\bigr\}}{1-\hat e_{-k}(X_i)}
+\hat Q^{0}_{\tau,-k}(X_i,q)
\right] \\
&\quad
+\frac{(1-T_i)(1-G_i)}{\nu}\hat\alpha_{0,-k}(S_i,X_i)
\left\{
\rho_\tau(Y_i-q)-\hat Q^{0}_{\tau,-k}(S_i,X_i,q)
\right\}.
\end{align*}
The doubly robust estimator of the quantile treatment effects is then
\[
\hat\Delta_{dr}
=
\hat q_{dr}^1-\hat q_{dr}^0.
\]
\end{definition}

\iffalse
As defined in subsection~\ref{sec:Semiparametric Efficiency}, $\psi_t(W_i;q,\hat\eta_{-k})$ is the derivative of $\Psi_t(W_i;q,\hat\eta_{-k})$ with respect to $q$. Therefore, solving the optimization problem in Definition~\ref{def:dr estimator} is equivalent to solving the estimating equation

\begin{align*}
\frac{1}{N}\sum_{k=1}^K\sum_{i\in\mathcal I_k}\psi_t(W_i;q,\hat\eta_{-k})=0,
\end{align*}

whose solution is $\hat q_{dr}^t$.
\fi

\subsection{Asymptotic properties}

In this subsection, we establish the consistency, asymptotic normality, and semiparametric efficiency of the proposed DR estimator. To this end, We first introduce some regularity conditions required to establish the asymptotic properties of our proposed DR estimator. 

Let $\| \cdot \|_{P,q}$ denote the $L^q(P)$ norm$\|f(W)\|_{P,q} = \left( \int |f(\omega)|^q \, dP(\omega) \right)^{1/q}$ and Let $\delta(F_1, F_2)$ be the total variation distance between two probability measures with CDFs $F_1$ and $F_2$.

\begin{assumption}\label{ass:regularity}
For t=0,1,

   \begin{itemize}
    \item[\textbf{(A1)}]\label{A1}
    $\left\| \widehat{e}(X) - e(X) \right\|_{P,2} = o(n^{-1/4})$ and $\left\| \widehat{g}_t(S,X) - g_t(S,X) \right\|_{P,2} = o(n^{-1/4}).$

    \item[\textbf{(A2)}]\label{A2}
    $\left\| \delta \left(\hat{F}_t(q\mid S,X),F_t(q\mid S,X) \right) \right\|_{P,2} = o(n^{-1/4}).$
    
    \item[\textbf{(A3)}]\label{A3}
    The conditional density $f_t(q\mid S,X)$ is uniformly bounded.
    
    \item[\textbf{(A4)}]\label{A4}
    There exist a constant $C > 0$, such that $f_t(q\mid S,X)\geq C$ holds for all $q$ in a neighborhood of $q^t_{\tau}$.
    
    \item[\textbf{(A5)}]\label{A5}
     $\exists$ $\varepsilon > 0$, such that $\Pr \left( \varepsilon \leq \hat{e}(X) \leq 1 - \varepsilon \right) = 1$ and $\Pr \left( \varepsilon \leq \hat{g}_t(S,X) \leq 1 - \varepsilon \right) = 1.$
    
    \item[\textbf{(A6)}]\label{A^}
    The derivative \(\partial_q f_t(q\mid S,X) \big|_{q = q^t_\tau}\) is bounded.
\end{itemize} 
\end{assumption}

Assumptions (A1)-(A2) impose rate conditions on the estimation errors of the propensity score and the conditional outcome distribution, respectively. Assumptions (A3)-(A4) provide regularity conditions on the conditional density: (A3) requires the conditional pdfs to be uniformly bounded, while (A4) ensures the identifiability of the quantile $q^t_\tau$, thereby guaranteeing the well-posedness of the quantile estimation problem and excluding degenerate cases. Assumption (A5) enforces the overlap condition for the estimated propensity score, and (A6) controls the smoothness of the density in a neighborhood of the true quantile $q^t_{\tau}$.

The following result establishes the doubly robust property of our estomator \(\hat\Delta_{dr}\).

\begin{lemma}\label{dr property}
    Under Assumptions~\ref{unconfoundedness}--\ref{Monotonicity} and Assumptions~(A4)--(A5) in Assumption~\ref{ass:regularity}, $\hat{\Delta}_{dr}$ is a consistent estimator of $\Delta_\tau$ as long as one of following two parts of models is consistently estimated:
    \begin{enumerate}
        \item The score functions(propensity score and sample score), i.e. $e(X)$ and $g_t(S,X)$.
        \item The conditional expectation functions, $F_t(q^t_\tau\mid S,X)$ and $m_t(X,q^t_\tau)$.
    \end{enumerate}
\end{lemma}

Thus, consistency of the proposed estimator does not require all nuisance functions to be consistently estimated, nor does it require prior knowledge of which components are correctly specified. In this sense, the estimator is robust to misspecification of some nuisance functions, provided that the remaining components are estimated consistently.

We now establish the asymptotic normality of the proposed estimator $\hat{\Delta}_{dr}$.

\begin{theorem}[Asymptotic properties of $\hat{\Delta}_{dr}$]\label{asymtotic normality}
Let \(\phi_i\) be the sample value of \(\phi\) at the \(i\)-th observation, for \(i=1,\dots,N\). Under Assumptions in Lemma~\ref{dr property} and the rest part of Assumption~\ref{ass:regularity}, the estimator $\hat{\Delta}_{dr}$ satisfies
\begin{align*}
    \sqrt{N}(\hat{\Delta}_{dr} - \Delta_\tau) 
    = -\frac{1}{\sqrt{N}} \sum_{i=1}^{N}\phi_i + o_p(1) 
    \xrightarrow{d} \mathcal{N}(0, V),
\end{align*}
where $V=\E[\phi^2]$  coincides with the semiparametric efficiency bound of $\Delta_\tau$.
\end{theorem}

Theorem~\ref{asymtotic normality} establishes that the proposed estimator $\hat{\Delta}_{dr}$ is $\sqrt{N}$-consistent and asymptotically normal, and that it achieves the semiparametric efficiency bound. This implies that no regular estimator can asymptotically improve upon its variance.

Importantly, these results hold under relatively mild conditions on the estimation of the nuisance functions. Specifically, all nuisance estimators are only required to converge at rate $o_p(N^{-1/4})$ in $L^2(P)$ norm—significantly slower than the standard parametric rate $O_p(N^{-1/2})$. Moreover, thanks to cross-fitting, we impose no Donsker class or uniform entropy conditions, thereby allowing the use of modern, flexible machine learning or deep learning methods that may not satisfy classical empirical process constraints\parencite{chernozhukov2018double}.

In addition, the asymptotic normality result relies on a product-rate condition, meaning that the product of estimation errors for pairs of nuisance functions must converge faster than $O_p(N^{-1/2})$. This allows slower convergence of some components to be compensated by faster convergence of others, enhancing robustness in practice and enabling greater flexibility in the choice of learners for different parts of the nuisance structure.

\subsection{Variance estimation}

In the previous subsection, we established the asymptotic normality of the DR estimator under mild conditions. This suggests that if we can estimate its asymptotic variance, then we can construct confidence intervals on $\hat\Delta_\tau$. In this subsection we provide a variance estimator and prove its consistency which immediately lends itself to confidence interval construction.

\begin{definition}[Variance estimator]
For each \(i=1,\dots,N\), define
\[
\hat\phi_i
=
\hat{J}_1^{-1}\psi\!\left(W_i;\hat q^1_{dr},\hat\eta_{k}\right)
-
\hat{J}_0^{-1}\psi\!\left(W_i;\hat q^0_{dr},\hat\eta_{k}\right),
\]
where \(\hat J_t\) is the estimator of $J_t$ corresponding to \(q_\tau^t\). Then define
\[
\hat V
=
\frac{1}{N}\sum_{i=1}^N
\left(\hat\phi_i-\bar{\hat\phi}\right)^2,
\qquad
\bar{\hat\phi}
=
\frac{1}{N}\sum_{i=1}^N \hat\phi_i .
\]
\end{definition}

We next establish the consistency of $\hat{V}$ and then establish the confidence interval which relies on another assumption. 

\begin{assumption}\label{ass: consistency of J}
 For \(t=0,1\), \(
\left\|
\hat f_t(q_\tau^t\mid S,X)-f_t(q_\tau^t\mid S,X)
\right\|_{P,2}
=o_p(1).
\) 
\end{assumption}

This assumption imposes \(L^2(P)\)-consistency of the density estimator \(\hat f_t\) at the target quantile \(q_\tau^t\). As a consequence, it justifies the use of a consistent estimator \(\hat J_t\) for \(J_t\). However, verifying this condition can be challenging in practice, since it requires estimating \(f_t\), which is a density estimation problem. We discuss our estimation strategy in detail in Remark~\ref{rem: estimating detail}.

\begin{theorem}[Confidence Interval]\label{Confidence Interval}
Under the assumptions in Lemma~\ref{dr property} and Assumption~\ref{ass: consistency of J}, the variance estimator
\[
\hat{V}\xrightarrow{p} V \quad \text{as } N \to \infty,
\]

is consistent for the asymptotic variance $V$.
Consequently, the following $(1 - \alpha) \times 100\%$ confidence interval
\begin{align*}
\mathrm{CI} = \left( \hat{\Delta}_\tau - z_{1-\alpha/2}\cdot \sqrt{\hat{V}/N},\;
                  \hat{\Delta} + z_{1-\alpha/2}\cdot \sqrt{\hat{V}/N} \right)
\end{align*}
satisfies
\[
\mathbb{P}(\Delta_\tau \in \mathrm{CI}) \to 1 - \alpha \quad \text{as } N \to \infty,
\]
where $z_{1-\alpha/2}$ is the $(1 - \alpha/2)$-quantile of the standard normal distribution.
\end{theorem}

Theorem~\ref{Confidence Interval} demonstrates that, under the same set of assumptions, the semiparametric efficiency bound $\mathbb{E}[\phi^2]$ can be consistently estimated by constructing its empirical analogue using cross-fitted nuisance estimators and the proposed estimator $\hat{\Delta}_{dr}$. This result ensures that valid inference for $\Delta_\tau$ can be conducted without requiring knowledge of the true nuisance functions or imposing stringent structural assumptions.

\begin{remark}\label{rem: estimating detail}
In practice, we can estimate the nuisance functions using flexible machine learning or deep learning methods. 
Specifically, the propensity score \(e(X)\) and the sample score \(g_t(S,X)\) can be estimated using any supervised learning method. In our simulation study, We employ generalized random forest(GRF) for theses components\parencite{athey2019generalized}.
For outcome-related nuisance functions, one may in principle apply supervised learning methods separately for each value of \(q\). However, such a method is computationally intensive when the search space of \(q\) is large. To improve computational efficiency, we instead fit a mixture density network (MDN) on the observational sample to estimate the conditional distribution of \(Y\mid S,X,T=t,G=0\). This approach directly yields estimators of \(F_t(q\mid S,X)\), \(f_t(q\mid S,X)\), and \(Q_\tau^t(S,X,q)\). To approximate the transported quantities in the experimental sample, we additionally fit an MDN for the surrogate distribution \(S\mid X,T=t,G=1\). The transported quantities \(Q_\tau^t(X,q)\), \(m_t(X,q)\), and \(\mathbb{E}\!\left[f_t(q\mid S,X)\mid X,T=t,G=1\right]\) are then evaluated by nested Monte Carlo integration. The outer Monte Carlo step draws samples from the fitted distribution of \(S\mid X,T=t,G=1\), and the inner quantities are computed using the fitted model for \(Y\mid S,X,T=t,G=0\). Based on the resulting estimate of the transported density, \(J_t\) is finally estimated by plug-in averaging over the experimental sample.

\end{remark}

\section{Simulation}\label{sec: Simulation}

We conduct a simulation study to investigate the finite-sample performance of the proposed doubly robust (DR) estimator. In each simulation setting, we generate an RCT sample and an observational sample under a common structural model.

In the RCT sample, let \(X=(X_1,X_2)^\top \sim \mathcal N(0,I_2)\), and assign treatment by \(T\sim \mathrm{Bernoulli}(0.5)\) independently of \(X\). The surrogate is generated as
\[
S=2(X_1+X_2)+T+\varepsilon_S,
\qquad \varepsilon_S\sim \mathcal N(0,1),
\]
and the long-term outcome is generated as
\[
Y=T+3(X_1+X_2)+S+\varepsilon_Y.
\]

In the observational sample, covariates are drawn from
\[
X\sim \mathcal N(0.5\mathbf 1_2,\,1.5^2 I_2),
\]
and treatment is assigned according to
\[
\Pr(T=1\mid X,G=0)=\mathrm{expit}\{0.25X_1+0.25X_2\}.
\]
Conditional on \((X,T)\), both \(S\) and \(Y\) follow the same structural equations as in the RCT sample.

To examine performance under different tail behaviors, we consider
\[
\varepsilon_Y=\sigma_{\kappa}\, t(\kappa),
\qquad \kappa\in\{3,5,7,9\},
\]
where
\[
\sigma_{\kappa}=\sqrt{\frac{\kappa-2}{\kappa}},
\]
so that \(\varepsilon_Y\) has unit variance, and \(t(\kappa)\) denotes the \(t\)-distribution with \(\kappa\) degrees of freedom. We also consider the Gaussian case
\[
\varepsilon_Y\sim \mathcal N(0,1).
\]
We consider RCT sample sizes \(n_{\mathrm{rct}}\in\{500,1000,2000\}\), and set the observational sample size to \(n_{\mathrm{obs}}=5n_{\mathrm{rct}}\). For each configuration, we use \(K=5\) folds for cross-fitting and repeat the experiment for \(1000\) Monte Carlo replications. The target quantile levels are \(\tau\in\{0.25,0.50,0.75\}\).

The true quantile treatment effects used to evaluate Monte Carlo bias are approximated by \(10{,}000{,}000\) Monte Carlo draws under the same data-generating mechanism.

Each simulation study is based on $1000$ replications.
In the following table, Bias and SD are the Monte Carlo bias and standard deviation of the point estimates.
ESE and CP95 are the averages of the estimated asymptotic standard errors and the empirical coverage proportions of the nominal $95\%$ Wald confidence intervals based on influence function, respectively.
Table~\ref{tab:sim_main} reports Bias(SD), ESE and CP95 for $\tau\in\{0.25,0.50,0.75\}$ across tail-heaviness levels and RCT sample sizes.

\begin{table}[!htbp]
\centering
\small
\renewcommand{\arraystretch}{1.15}
\setlength{\tabcolsep}{6pt}
\caption{Finite-sample performance across tail levels and sample sizes. }
\begin{tabular}{llccc ccc ccc}
\toprule
\multirow{2}{*}{$\varepsilon_Y$} & \multirow{2}{*}{$n_{rct}$}
& \multicolumn{3}{c}{$\tau=0.25$}
& \multicolumn{3}{c}{$\tau=0.50$}
& \multicolumn{3}{c}{$\tau=0.75$} \\
\cmidrule(lr){3-5}\cmidrule(lr){6-8}\cmidrule(lr){9-11}
& & Bias (SD) & ESE & CP95 & Bias (SD) & ESE & CP95 & Bias (SD) & ESE & CP95 \\
\midrule

\multirow{3}{*}{$\sigma_3t(3)$}
& 500  & \text{0.009(0.93)} & \text{1.01} & \text{96.9} & \text{0.010(0.86)} & \text{0.90} & \text{95.8} & \text{0.002(0.91)} & \text{0.99} & \text{97.4} \\
& 1000 & \text{0.012(0.66)} & \text{0.72} & \text{96.8} & \text{0.020(0.61)} & \text{0.64} & \text{96.4} & \text{0.020(0.66)} & \text{0.70} & \text{96.8} \\
& 2000 & \text{0.018(0.47)} & \text{0.51} & \text{96.7} & \text{0.007(0.41)} & \text{0.46} & \text{96.8} & \text{0.004(0.45)} & \text{0.49} & \text{97.2} \\

\addlinespace[2pt]

\multirow{3}{*}{$\sigma_5t(5)$}
& 500  & \text{0.012(0.93)} & \text{1.01} & \text{96.8} & \text{0.010(0.86)} & \text{0.90} & \text{96.8} & \text{0.017(0.90)} & \text{0.99} & \text{96.9} \\
& 1000 & \text{0.002(0.66)} & \text{0.71} & \text{96.5} & \text{0.019(0.61)} & \text{0.64} & \text{96.0} & \text{0.008(0.65)} & \text{0.70} & \text{95.5} \\
& 2000 & \text{0.009(0.46)} & \text{0.51} & \text{96.4} & \text{0.009(0.40)} & \text{0.45} & \text{97.6} & \text{0.006(0.44)} & \text{0.49} & \text{96.7} \\

\addlinespace[2pt]

\multirow{3}{*}{$\sigma_7t(7)$}
& 500  & \text{0.042(0.93)} & \text{1.02} & \text{97.5} & \text{0.020(0.86)} & \text{0.90} & \text{96.1} & \text{0.001(0.92)} & \text{0.98} & \text{97.0} \\
& 1000 & \text{0.001(0.66)} & \text{0.72} & \text{97.2} & \text{0.018(0.60)} & \text{0.64} & \text{95.8} & \text{0.004(0.65)} & \text{0.70} & \text{96.0} \\
& 2000 & \text{0.008(0.46)} & \text{0.51} & \text{96.8} & \text{0.002(0.40)} & \text{0.45} & \text{97.5} & \text{0.001(0.44)} & \text{0.49} & \text{97.4} \\

\addlinespace[2pt]

\multirow{3}{*}{$\sigma_9t(9)$}
& 500  & \text{0.004(0.94)} & \text{1.02} & \text{97.6} & \text{0.017(0.87)} & \text{0.90} & \text{95.2} & \text{0.006(0.92)} & \text{0.98} & \text{96.9} \\
& 1000 & \text{0.002(0.67)} & \text{0.72} & \text{96.5} & \text{0.011(0.59)} & \text{0.64} & \text{96.7} & \text{0.007(0.65)} & \text{0.70} & \text{96.1} \\
& 2000 & \text{0.013(0.46)} & \text{0.51} & \text{96.4} & \text{0.008(0.40)} & \text{0.45} & \text{97.4} & \text{0.001(0.44)} & \text{0.49} & \text{96.4} \\

\addlinespace[2pt]

\multirow{3}{*}{$\mathcal N(0,1)$}
& 500  & $0.017 $ ($0.93$) & $1.02$ & $97.1$ & \text{0.000(0.86)} & \text{0.91} & \text{96.3} & \text{0.012(0.92)} & \text{0.98} & \text{97.0} \\
& 1000 & $0.006 $ ($0.65$) & $0.72$ & $96.7$ & \text{0.021(0.60)} & \text{0.64} & \text{96.5} & \text{0.007(0.63)} & \text{0.70} & \text{97.2} \\
& 2000 & $0.008$ ($0.46$) & $0.51$ & $96.8$ & \text{0.002(0.40)} & \text{0.45} & \text{97.5} & \text{0.001(0.44)} & \text{0.49} & \text{97.4} \\
\addlinespace[2pt]
\bottomrule
\end{tabular}
\label{tab:sim_main}
\end{table}

Table~\ref{tab:sim_main} shows that the proposed estimator performs stably across all simulation settings. The empirical bias is negligible relative to the Monte Carlo standard deviation in every scenario, and the variability decreases substantially as the RCT sample size increases. The EIF-based standard error tracks the empirical standard deviation well, although it is slightly conservative, which leads to coverage probabilities that are consistently close to, and often slightly above, the nominal 95\% level. This mild conservativeness is likely driven by downward bias in the estimated density at the target quantile, which enters the EIF-based variance formula in the denominator. underestimating this quantity inflates the estimated variance and hence the standard error. Such behavior is common in quantile treatment effect inference, where variance estimation is well known to be sensitive to the estimation of the density at the target quantile; see, for example, \textcite{kallus2024localized}. Moreover, the finite-sample performance is largely insensitive to the tail heaviness of the outcome error distribution over the range of $t$-distributions and the Gaussian case considered here.

\section{Real-data application}\label{sec:Application}

IgA nephropathy (IgAN), also known as Berger's disease, is one of the most common primary chronic glomerular diseases worldwide \parencite{haas1997histologic}. Renin--angiotensin--aldosterone system (RAAS) inhibition is a standard treatment for IgAN and mainly delays disease progression by reducing proteinuria and blood pressure \parencite{zhang2021effects}. However, even with RAAS therapy, patients with IgAN remain at risk of kidney failure \parencite{liu2019effects}. As an immunomodulatory agent, hydroxychloroquine (HCQ) has recently been considered a potential treatment option for IgAN. Existing studies have shown that, compared with RAAS therapy alone, HCQ combined with RAAS therapy can more effectively reduce proteinuria within 6 months \parencite{yang2018effects}. Nevertheless, evidence on the effect of HCQ on long-term renal outcomes remains limited \parencite{zhang2021effects}. Therefore, it is of substantial clinical interest to further evaluate the effect of HCQ on long-term renal outcomes by combining randomized trial data with external observational data.

The experimental data come from a double-blind, randomized, placebo-controlled trial comparing HCQ plus optimized RAAS inhibition with standard RAAS therapy alone. In this paper, we combine the trial data with external observational data collected from the same hospital.

After data quality control, the final analytic sample consists of 52 patients from the experimental sample and 474 patients from the observational sample. In the experimental sample, 28 patients were assigned to the control group and 24 to the treatment group. In the observational sample, 418 patients were in the control group and 56 patients were in the treatment group.

The long-term outcome is defined as the log decline in glomerular filtration rate (GFR) from baseline to the minimum observed GFR during follow-up, namely
\[
\log(GFR_0)-\log(GFR_{\min}),
\]
so that a larger value indicates more severe deterioration in renal function. Because the trial lasted only 6 months, the long-term renal outcome is unavailable in the experimental sample. The surrogate is defined as the percentage change in proteinuria from baseline to 6 months. The baseline covariates are aligned across the experimental and observational data and include sex, age, baseline proteinuria, baseline GFR, and several pathological indicators associated with kidney failure \parencite{shi2011pathologic}.

\textcite{liu2019effects} showed that, compared with RAAS therapy alone, HCQ combined with RAAS therapy achieved a more favorable effect on the surrogate outcome. Furthermore, under the same clinical setting, \textcite{hu2025identification} found that HCQ combined with RAAS therapy had a more favorable average treatment effect on the long-term renal outcome than RAAS therapy alone. Building on these findings, we further investigate the quantile treatment effect on the long-term outcome in the target trial population, in order to assess whether the new treatment regimen performs better than standard therapy at different parts of the long-term outcome distribution.

\begin{figure}[ht]
\centering
\includegraphics[width=0.8\textwidth]{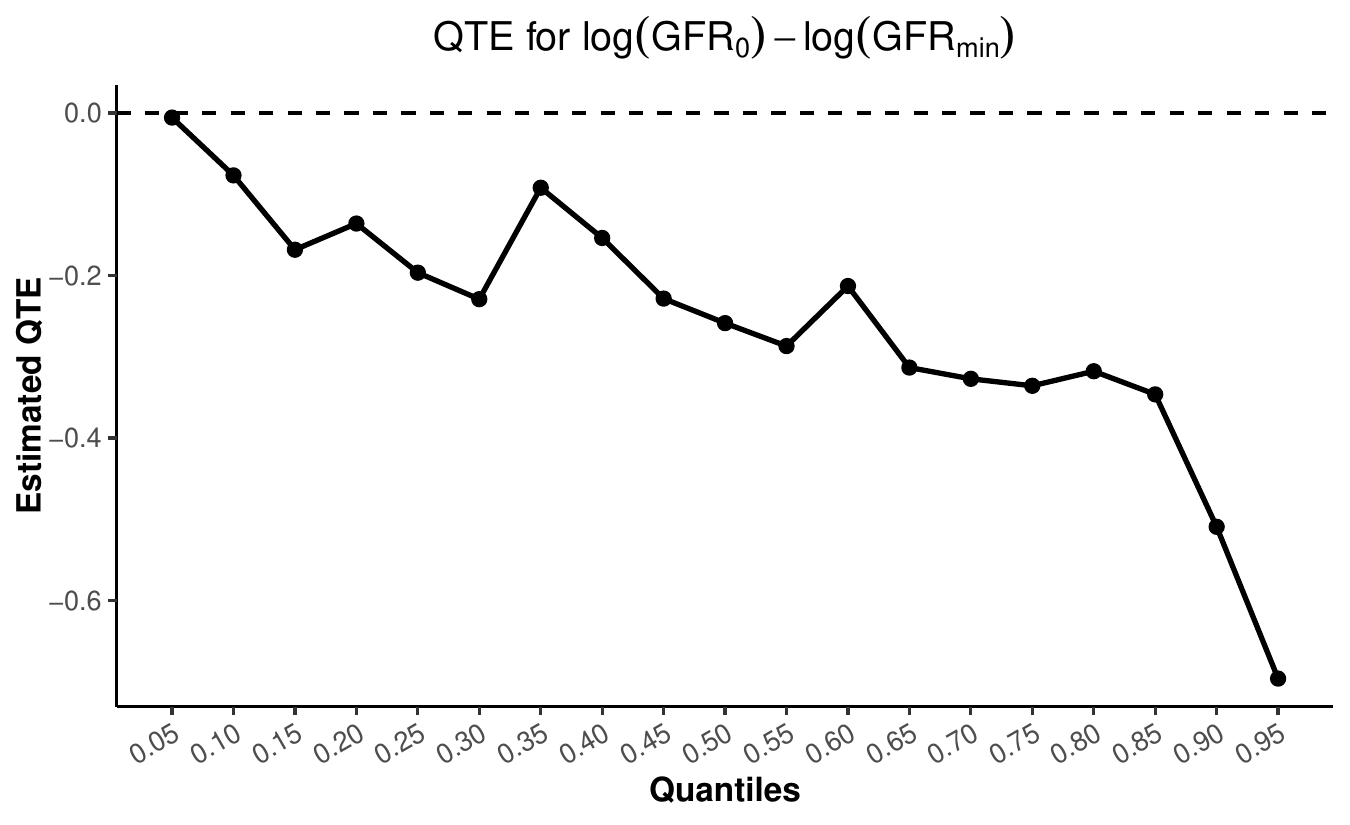}
\caption{
Estimated quantile treatment effect curve for \(\log(GFR_0)-\log(GFR_{\min})\). The solid black curve represents the estimated QTE. The horizontal dashed line corresponds to zero treatment effect.
}
\label{fig:application_qte}
\end{figure}

Figure~\ref{fig:application_qte} presents the estimated QTE curve over quantile levels from \(\tau=0.05\) to \(\tau=0.95\). The estimated QTE remains below zero over most of the quantile range considered. This suggests that, in the target trial population, HCQ combined with RAAS therapy is associated with a lower degree of long-term renal deterioration relative to RAAS therapy alone.

Moreover, the magnitude of the estimated treatment effect becomes larger at higher quantiles. Specifically, the estimated effect is relatively close to zero in the lower part of the outcome distribution, becomes more pronounced around the middle quantiles, and is strongest in the upper tail. This pattern suggests that the treatment may provide greater protection for patients who would otherwise experience more severe renal decline.

Overall, this application further suggests, beyond previous studies, that HCQ combined with RAAS therapy may not only provide renal protection on average, but also exhibit substantial heterogeneity across different parts of the outcome distribution. In particular, for patients who would otherwise experience more severe renal deterioration, the treatment effect may be more pronounced.

\section{Conclusion}\label{sec:Conclusion}
We study the estimation of quantile treatment effects for long-term outcomes when the primary outcome is unavailable in the randomized trial but observed in an external observational sample. Our main contribution is to extend surrogate-based data integration from average treatment effects to quantile treatment effects. To this end, we derive the efficient influence function and develop a doubly robust estimator that allows flexible estimation of nuisance functions.

Under standard identification conditions, the target quantiles are identifiable. The proposed estimator remains consistent if either the score-related nuisance functions or the outcome-regression-related nuisance functions are consistently estimated, and it supports valid large-sample inference under regularity conditions. These results provide a practical framework for estimating long-term quantile treatment effects when the primary outcome is not observed in randomized trials.

In practice, quantile inference depends on estimating the density around the target quantile, which may be unstable in finite samples, especially in the tails. Future work may consider settings with partial long-term outcome observation in the trial, longitudinal surrogate information, and sensitivity analysis for departures from the identifying conditions.
\printbibliography

@misc{chernozhukov2018double,
  title={Double/debiased machine learning for treatment and structural parameters},
  author={Chernozhukov, Victor and Chetverikov, Denis and Demirer, Mert and Duflo, Esther and Hansen, Christian and Newey, Whitney and Robins, James},
  year={2018},
  publisher={Oxford University Press Oxford, UK}
}

@article{athey2025surrogate,
  title={The surrogate index: Combining short-term proxies to estimate long-term treatment effects more rapidly and precisely},
  author={Athey, Susan and Chetty, Raj and Imbens, Guido W and Kang, Hyunseung},
  journal={Review of Economic Studies},
  pages={rdaf087},
  year={2025},
  publisher={Oxford University Press UK}
}

@article{imbens2025long,
  title={Long-term causal inference under persistent confounding via data combination},
  author={Imbens, Guido and Kallus, Nathan and Mao, Xiaojie and Wang, Yuhao},
  journal={Journal of the Royal Statistical Society Series B: Statistical Methodology},
  volume={87},
  number={2},
  pages={362--388},
  year={2025},
  publisher={Oxford University Press UK}
}

@article{kallus2025role,
  title={On the role of surrogates in the efficient estimation of treatment effects with limited outcome data},
  author={Kallus, Nathan and Mao, Xiaojie},
  journal={Journal of the Royal Statistical Society Series B: Statistical Methodology},
  volume={87},
  number={2},
  pages={480--509},
  year={2025},
  publisher={Oxford University Press UK}
}

@article{chen2023semiparametric,
  title={Semiparametric estimation of long-term treatment effects},
  author={Chen, Jiafeng and Ritzwoller, David M},
  journal={Journal of Econometrics},
  volume={237},
  number={2},
  pages={105545},
  year={2023},
  publisher={Elsevier}
}

@article{kallus2024localized,
  title={Localized debiased machine learning: Efficient inference on quantile treatment effects and beyond},
  author={Kallus, Nathan and Mao, Xiaojie and Uehara, Masatoshi},
  journal={Journal of Machine Learning Research},
  volume={25},
  number={16},
  pages={1--59},
  year={2024}
}

@book{HernanRobins2020,
  author    = {Hern{\'a}n, Miguel A. and Robins, James M.},
  title     = {Causal Inference: What If},
  year      = {2020},
  address   = {Boca Raton},
  publisher = {Chapman \& Hall/CRC},
  url       = {https://www.hsph.harvard.edu/miguel-hernan/causal-inference-book/},
  note      = {Free eBook available online}
}

@article{rosenbaum1983central,
  title={The central role of the propensity score in observational studies for causal effects},
  author={Rosenbaum, Paul R and Rubin, Donald B},
  journal={Biometrika},
  volume={70},
  number={1},
  pages={41--55},
  year={1983},
  publisher={Oxford University Press}
}

@article{firpo2007efficient,
  title={Efficient semiparametric estimation of quantile treatment effects},
  author={Firpo, Sergio},
  journal={Econometrica},
  volume={75},
  number={1},
  pages={259--276},
  year={2007},
  publisher={Wiley Online Library}
}

@article{frolich2013unconditional,
  title={Unconditional quantile treatment effects under endogeneity},
  author={Fr{\"o}lich, Markus and Melly, Blaise},
  journal={Journal of business \& economic statistics},
  volume={31},
  number={3},
  pages={346--357},
  year={2013},
  publisher={Taylor \& Francis}
}

@article{diaz2017efficient,
  title={Efficient estimation of quantiles in missing data models},
  author={D{\'\i}az, Iv{\'a}n},
  journal={Journal of Statistical Planning and Inference},
  volume={190},
  pages={39--51},
  year={2017},
  publisher={Elsevier}
}

@article{cheng2024doubly,
  title={Doubly robust estimation and sensitivity analysis for marginal structural quantile models},
  author={Cheng, Chao and Hu, Liangyuan and Li, Fan},
  journal={Biometrics},
  volume={80},
  number={2},
  pages={ujae045},
  year={2024},
  publisher={Oxford University Press}
}

@article{hines2022demystifying,
  title={Demystifying statistical learning based on efficient influence functions},
  author={Hines, Oliver and Dukes, Oliver and Diaz-Ordaz, Karla and Vansteelandt, Stijn},
  journal={The American Statistician},
  volume={76},
  number={3},
  pages={292--304},
  year={2022},
  publisher={Taylor \& Francis}
}

@article{hu2025identification,
  title={IDENTIFICATION AND ESTIMATION OF TREATMENT EFFECTS ON LONG-TERM OUTCOMES IN CLINICAL TRIALS WITH EXTERNAL OBSERVATIONAL DATA},
  author={Hu, Wenjie and Zhou, Xiao-Hua and Wu, Peng},
  journal={Statistica Sinica},
  volume={35},
  pages={959--980},
  year={2025}
}

@article{haas1997histologic,
  title={Histologic subclassification of IgA nephropathy: a clinicopathologic study of 244 cases},
  author={Haas, Mark},
  journal={American journal of kidney diseases},
  volume={29},
  number={6},
  pages={829--842},
  year={1997},
  publisher={Elsevier}
}

@article{liu2019effects,
  title={Effects of hydroxychloroquine on proteinuria in IgA nephropathy: a randomized controlled trial},
  author={Liu, Li-Jun and Yang, Ya-zi and Shi, Su-Fang and Bao, Yun-Fei and Yang, Chao and Zhu, Sai-Nan and Sui, Gui-Li and Chen, Yu-Qing and Lv, Ji-Cheng and Zhang, Hong},
  journal={American Journal of Kidney Diseases},
  volume={74},
  number={1},
  pages={15--22},
  year={2019},
  publisher={Elsevier}
}

@article{zhang2021effects,
  title={Effects of Hydroxychloroquine on Proteinuria in IgA Nephropathy: A Systematic Review and Meta-Analysis},
  author={Zhang, Jialing and Lu, Xiangxue and Feng, Jianan and Li, Han and Wang, Shixiang},
  journal={BioMed Research International},
  volume={2021},
  number={1},
  pages={9171715},
  year={2021},
  publisher={Wiley Online Library}
}

@article{yang2018effects,
  title={Effects of hydroxychloroquine on proteinuria in immunoglobulin A nephropathy},
  author={Yang, Ya-Zi and Liu, Li-Jun and Shi, Su-Fang and Wang, Jin-Wei and Chen, Yu-Qing and Lv, Ji-Cheng and Zhang, Hong},
  journal={American journal of nephrology},
  volume={47},
  number={3},
  pages={145--152},
  year={2018},
  publisher={S. Karger AG Basel, Switzerland}
}

@article{shi2011pathologic,
  title={Pathologic predictors of renal outcome and therapeutic efficacy in IgA nephropathy: validation of the oxford classification},
  author={Shi, Su-Fang and Wang, Su-Xia and Jiang, Lei and Lv, Ji-Cheng and Liu, Li-Jun and Chen, Yu-Qing and Zhu, Sai-Nan and Liu, Gang and Zou, Wan-Zhong and Zhang, Hong and others},
  journal={Clinical Journal of the American Society of Nephrology},
  volume={6},
  number={9},
  pages={2175--2184},
  year={2011},
  publisher={LWW}
}

@article{rubin1974estimating,
  title={Estimating causal effects of treatments in randomized and nonrandomized studies.},
  author={Rubin, Donald B},
  journal={Journal of educational Psychology},
  volume={66},
  number={5},
  pages={688},
  year={1974},
  publisher={American Psychological Association}
}

@article{splawa1990application,
  title={On the application of probability theory to agricultural experiments. Essay on principles. Section 9},
  author={Splawa-Neyman, Jerzy and Dabrowska, Dorota M and Speed, Terrence P},
  journal={Statistical Science},
  pages={465--472},
  year={1990},
  publisher={JSTOR}
}

@article{athey2019generalized,
  title={Generalized random forests},
  author={Athey, Susan and Tibshirani, Julie and Wager, Stefan},
  journal={The Annals of Statistics},
  volume={47},
  number={2},
  pages={1148--1178},
  year={2019}
}

\appendix

%\begin{refsection}
%\subfile{appendix}

%\newrefcontext[labelprefix={}]
%\printbibliography[heading=subbibliography,resetnumbers=true]
%\end{refsection}

\end{document}